\documentclass[11pt,twoside, final]{amsart}
\copyrightinfo{0}{Iranian Mathematical Society}
\usepackage{amsmath,amsthm,amscd,amsfonts,amssymb,enumerate}
\usepackage{graphicx}		
\usepackage{color}
\usepackage[colorlinks]{hyperref}

\theoremstyle{definition}

\theoremstyle{remark}

\numberwithin{equation}{section}
 
 \commby{Hamid Pezeshk}
\begin{document}

 
\title[ICSRBF FOR VOLTERRA'S POPULATION MODEL]{Collocation Method using Compactly Supported Radial Basis Function for Solving Volterra's Population Model} 
 
\author[Parand]{Kourosh Parand $^*$}
\address[First Author]{Department of Computer Sciences, Shahid Beheshti University, G.C. Tehran 19697-64166, Iran}
\email{k\_parand@sbu.ac.ir}

\author[Hemami]{Mohammad Hemami}
\address[Second Author]{Department of Computer Sciences, Shahid Beheshti University, G.C. Tehran 19697-64166, Iran}
\email{mohammadhemami@yahoo.com}

  \thanks{$^*$Corresponding author:Member of research group of Scientific Computing. Fax:
+98 2122431653.}
%
\maketitle
%

\begin{abstract}
In this paper, indirect collocation approach based on compactly supported radial basis function (CSRBF) is applied for solving Volterra's population model. The method reduces the solution of this problem to the solution of a system of algebraic equations. Volterra's model is a non-linear integro-differential equation where the integral term represents the effect of toxin. To solve the problem, we use the well-known CSRBF: $Wendland_{3,5}$. Numerical results and residual norm ($\|R(t)\|^2$) show good accuracy and rate of convergence. 

\end{abstract}
 
 \section{Introduction}
\label{intro}
The volterra's model for population growth of a species within a closed system is given in \cite{Scudo,TeBeest} as 
\begin{equation}
\frac{dp}{dt}=ap- bp^2 - cp \int_0^t p(x)dx,~~~~p(0)=p_0,
\label{Eq1}
\end{equation}
where $a>0$ is the birth rate coefficient, $b>0$ is the crowding coefficient and $c>0$ is the toxicity coefficient. the coefficient $c$ indicates the essential behaviour of the population evolution before its level falls to zero in the long term.  $p_0$ is the initial population and $p=p(t)$ denotes the population at time $t$. This model includes the well-known terms of a logistic equation, and in addition it, includes an integral term $cp \int_0^t p(x) dx$ that characterizes the accumulated toxicity produced since time zero\cite{TeBeest,Wazwaz}. \\
We apply scale time and population by introducing the non-dimensional variables 
\begin{equation}
t=\frac{tc}{b},~~~u=\frac{pb}{a},
\label{Eq2}
\end{equation}
to obtain the non-dimensional problem
\begin{equation}
\kappa\frac{du}{dt}=u-u^2-u\int_0^t u(x)dx,~~~~u(0)=u_0,
\label{Eq3}
\end{equation}
where $u(t)$ is the scaled population of identical individuals at time $t$ and $\kappa=\frac{c}{ab}$ is a prescribed non-dimensional parameter. The only equilibrium solution of Eq. (\ref{Eq3}) is the trivial solution $u(t)=0$ and the analytical solution \cite{TeBeest}
\begin{equation}
u(t)=u_0\exp(\frac{1}{\kappa}\int_0^t[1-u(\tau)-\int_0^\tau u(x)dx]d\tau)
\label{Eq4}
\end{equation}
shows that $u(t)>0$ for all $t$ if $u_0>0$.\par 
The solution of Eq. (\ref{Eq1}) has been of considerable concern. Although a closed form solution has been achieved in \cite{Scudo,Small}, it was formally shown that the closed form solution cannot lead to any insight into the behaviour of the population evolution \cite{Scudo}. Some approximate and numerical solutions for Volterra's Population Model have been reported. the successive approximations method was suggested for the solution of Eq. (\ref{Eq3}), but was not implemented. In this case, the solution $u(t)$ has a smaller amplitude compared to the amplitude of $u(t)$ for the case $\kappa \ll 1$.\par 
In \cite{TeBeest}, several numerical algorithms namely the Euler method, the modified Euler method, the classical fourth-order Runge-Kutta method and Runge-Kutta-Fehlberg method for the solution of Eq. (\ref{Eq3}) are obtained.  Moreover, a phase-plane analysis is implemented. In \cite{TeBeest}, the numerical results are correlated to give insight on the problem and its solution without using perturbation techniques. However, the performance of of  the traditional  numerical techniques is well-known in that it using provides grid points only, and in addition, it requires a large amounts of calculations.\par 
In \cite{Wazwaz}, the series solution method and the decomposition method are implemented independently  to Eq. (\ref{Eq3}) and to a related
non-linear ordinary differential equation. Furthermore, the Pad\'{e}approximations are used in the analysis to capture the essential behaviour of the populations $u(t)$ of identical individuals and approximation of $u_{max}$ and exact value of $u_{max}$ for different $\kappa$ were compared. 
Small \cite{Small} solved the Volterra’s Population Model by the singular perturbation method. This author scaled out the parameters of Eq. (\ref{Eq1}) as much as possible and considered four different ways to do this. He considered two cases $\kappa=\frac{c}{ab}$ small and $\kappa=\frac{c}{ab}$ large.\par
It is shown in \cite{Small} that for the case $\kappa \ll 1$, where populations are weakly sensitive to toxins, a rapid rise occurs along the logistic curve that will reach a peak and then is followed by a slow exponential decay. And, for large$\kappa$, where populations are
strongly sensitive to toxins, the solutions are proportional to $sech^2(t)$.\par 
In \cite{Al-Khaled} Adomian decomposition method and Sinc-Galerkin method were compared for the solution of the same integral
equation in this paper. This showed that Adomian decomposition method is more efficient and easier to use for the solution
of Volterra’s Population Model.\par 
In \cite{Ramezani} the approach is based upon composite spectral functions approximations. The properties of composite spectral
functions consisting of few terms of orthogonal functions utilized to reduce the solution of the Volterra’s model to the solution
of a system of algebraic equations.\par  
In \cite{Parand1} Rational Chebyshev and Hermite functions collocation approach were compared for the solution of Volterra’s
Population Model growth model of a species within a closed system. They reduced the solution of this problem to the solution
of a system of algebraic equations.\par 
In \cite{Parand6} applied two common collocation approaches based on radial basis functions to solve Volterra’s Population Model.\par 
In \cite{Marzban} a numerical method based on Hybrid function approximations was proposed to solve Volterra’s Population Model.
These Hybrid functions consist of block-pulse and Lagrange-interpolating polynomials.\par 
Also, in \cite{Parand7} Volterra’s population growth model of a species within a closed system is approximated by collocation
method based on two orthogonal functions, Sinc and Rational Legendre functions.
Momani et al. \cite{Momani} and Xu \cite{Xu} used a numerical and analytical algorithm for approximate solutions of a fractional population
growth model, respectively. The first algorithm is based on Adomian decomposition method (ADM) with Padé
approximants and the second algorithm is based on homotopy analysis method (HAM).

\section{ ICSRBF method}
\label{sec:4}
\subsection{CSRBF}
\label{ssec:41}
Many problems in science and engineering arise in infinite and semi-infinite domains. Different numerical methods have been proposed for solving problems on various domains such as FEM\cite{DH10,DH11}, FDM\cite{DH9,DH10} and Spectral\cite{DH5,Parand1} methods and meshfree method\cite{DH3,DH4}.
The use of the RBF is an one of the popular meshfree method for solving the differential equations  \cite{DH1,DH2}. For many years the global radial basis functions such as Gaussian, Multi quadric, Thin plate spline, Inverse multiqudric and etc was used  \cite{DH6,DH7,DH8}. These functions are globally supported and generate a system of equations  with ill-condition full matrix.To convert the ill-condition matrix to a well-condition matrix, CSRBFs can be used  instead of global RBFs. CSRBFs can convert the global scheme into a local one with banded matrices, Which makes the RBF method more feasible for solving large-scale problem \cite{Wong}.
\subsection*{Wendland's functions}
\label{ssec:42}
The most popular family of CSRBF are Wendland functions. This function introduced by Holger Wendland in 1995 \cite{Wendland}. Wendland starts with the truncated power function $\phi_l(r)=(1-r)^l_+$ which be strictly positive definite and radial on $\mathbb{R}^s$ for $l\geq \lfloor\frac{s}{2}\rfloor+1$ , and then he walks through dimension by repeatedly applying the operator I.\\
\textbf{Definition \cite{MeshfreeFasshhauer}}~~
with $\phi_l(r)=(1-r)^l_+$ we define
 \begin{equation}
 \phi_{s,k}=I^k\phi_{\lfloor\frac{s}{2}\rfloor +k+1},
 \end{equation}
 it turns out that the functions $\phi_{s,k}$ are all supported on [0,1].\\
\textbf{Theorem 1 \cite{MeshfreeFasshhauer}}~~
 The function $\phi_{s,k}$ are strictly positive definite (SPD)  and radial on $\mathbb{R}^s$ and are of the form
 \begin{equation}
\phi _{s,k}(r)=\begin{cases}
p_{s,k}(r)&   r\in [0,1],\nonumber\\
0&       r>1,
\end{cases}
\end{equation}
with a univariate polynomial $p_{s,k}$ of degree $\lfloor \frac{s}{2}\rfloor+3k+1$. Moreover, ،$\phi_{s,k}\in C^{2k}(R)$ are unique up to a constant factor, and the polynomial degree is minimal for given space dimension $s$ and smoothness $2k$ \cite{MeshfreeFasshhauer}.
Wendland gave recursive formulas for the functions $\phi_{s,k}$ for all $s, k$. We instead list the explicit formulas of \cite{Fasshauer}.\\
\textbf{Theorem 2 \cite{MeshfreeFasshhauer}}
~~~~The function $\phi_{s,k}$, $k=0,~1,~ 2,~ 3,$ have form
 \begin{eqnarray}
&&\phi_{s,0}=(1-r)^l_+,\nonumber\\
&&\phi_{s,1}\doteq(1-r) _+^{l+1}[(l+1)r+1],\nonumber\\
&&\phi_{s,2}\doteq(1-r)_+^{l+2}[(l^2+4l+3)r^2+(3l+6)r+3],\nonumber\\
&&\phi_{s,3}\doteq(1-r)_+^{l+3}[(l^3+9l^2+23l+15)r^3+(6l^2+36l+45)r^2\nonumber\\
&&+(15l+45)r+15],\nonumber
\end{eqnarray}
where $l=\lfloor \frac{s}{2}\rfloor +k+1$, and the symbol $\doteq$ denotes equality up to a multiplicative positive constant.\\
~~~~~~The case $k=0$ follows directly from the definition. application of the definition for the case $k=1$ yields
\begin{eqnarray}
&&\phi_{s,1}(r)=(I\phi_l)(r)\nonumber
=\int^\infty_r t\phi_l(t)dt  \nonumber\\
&&=\int^\infty _r t(1-t)^l_+dt\nonumber
=\int^1_r t(1-t)^ldt\nonumber\\
&&=\frac{1}{(l+1)(l+2)} (1-r)^{l+1}[(l+1)r+1],\nonumber
\end{eqnarray}
where the compact support of $\phi_l$ reduces the improper integral to a definite integral which can be evaluated using integration by parts. The other two cases are obtained similarly by repeated application of $I$.\cite{MeshfreeFasshhauer}
We showed the most of wendland functions in table \ref{Table. 1.} .
\begin{table}[htbpH]
\caption{\scriptsize
Wendland's compactly supported radial function for various choices of k and s=3.
}
\scriptsize \begin{tabular*}{\columnwidth}{@{\extracolsep{\fill}}*{3}{|c|}}
\hline
 $\phi_{s,k}$ & smoothness & SPD   \\
\hline
 $\phi_{3,0}(r) =(1-r)^2_+$ & $C^0$ & $\mathbb{R}^3$  \\
 \hline
 $\phi_{3,1}(r)\doteq(1-r)^4_+(4r+1)$ & $C^2$ & $\mathbb{R}^3$ \\
 \hline
 $\phi_{3,2}(r)\doteq(1-r)^6_+(35r^2+18r+3)$ & $C^4$ & $\mathbb{R}^3$  \\
 \hline
 $\phi_{3,3}(r)\doteq(1-r)^8_+(32r^3+25r^2+8r+1)$ & $C^6$ & $\mathbb{R}^3$  \\
 \hline
 $\phi_{3,4}(r)\doteq(1-r)^{10}_+(429r^4+450r^3+210r^2+50r+5)$ & $C^8$ & $\mathbb{R}^3$ \\
 \hline
 $\phi_{3,5}(r)\doteq(1-r)^{12}_+(2048r^5+2697r^4+1644r^3+566r^2+108r+9)$ & $C^{10}$ & $\mathbb{R}^3$  \\
 \hline
\end{tabular*}
\label{Table. 1.}
\end{table}

\subsection{Interpolation by CSRBFs}
\label{ssec:43}
The one-dimensional function $y(x)$ to be interpolated or approximated can be represented by an CSRBF as
\begin{equation}
y(x)\approx y_n(x)=\sum_{i=1}^N \xi_i \phi_i (x)=\Phi ^T(x)\Xi\nonumber,
\label{systmat}
\end{equation}
where
\begin{eqnarray}
&&\phi_i (x)=\phi (\dfrac{\|x-x_i\|}{r_\omega}),\nonumber\\
&&\Phi ^T(x)=[\phi_1 (x),\phi_2 (x),\cdots ,\phi _N(x)],\nonumber\\
&&\Xi =[\xi_1,\xi_2, \cdots , \xi_N ]^T,\nonumber
\end{eqnarray}
\begin{equation}
\end{equation}
$x$ is the input, $r_\omega$ is the local support domain and $\xi_i$s are the set of coefficients to be determined. By using the local support domain, we mapped the domain of problem to CSRBF local domain. By choosing $N$ interpolate nodes $(x_j,~ j=1, 2,\cdots,  N)$ in domain:
\begin{equation}
y_j=\sum_{i=1}^N\xi_i\phi_i (x_j)   (j=1, 2, \cdots, N).\nonumber
\end{equation}
To summarize the discussion on the coefficients matrix, we define
\begin{equation}
A\Xi = Y,
\label{system}
\end{equation}
where :
\begin{eqnarray}
&&Y=[y_1, y_2, \cdots,  y_N]^T,\nonumber\\
&&A=[\Phi ^T(x_1), \Phi ^T(x_2), \cdots , \Phi ^T(x_N)]^T\nonumber\\
&&=\begin{pmatrix} 
\phi_1(x_1)&\phi_2(x_1)&\cdots &\phi_N(x_1)\cr \phi_1(x_2)&\phi_2(x_2)&\cdots &\phi_N(x_2)\cr \vdots &\vdots &\ddots &\vdots \cr \phi_1(x_N)&\phi_2(x_N)&\cdots &\phi_N(x_N)\nonumber
\end{pmatrix}.
\end{eqnarray}
Note that $\phi_i(x_j)=\phi(\dfrac{\|x_i-x_j\|}{r_\omega})$, by solving the system (\ref{system}), the unknown coefficients $\xi_i$  will be achieved.
\subsection{ICSRBF method}
\label{ssec:44}
In the indirect method, the formulation of the problem starts with the decomposition of the highest order derivative under
consideration into CRBF. The obtained derivative expression is then integrated to yield expressions for lower order derivatives
and finally for the original function itself. \\
We approximate $\frac{du}{dt}$ for solving the model by ICSRBF:
\begin{equation}
\frac{du}{dt}\simeq \hat{u}_N(t)=\sum_{i=1}^N \xi_i \phi_i(t)= \Phi^T(t)\Xi,
\label{Eq5}
\end{equation}
by using integral operator $I_\vartheta f(t)=\int_0^t f(x)dx$ we have 
\begin{eqnarray}
\int_0^t du \simeq \int_0^t \hat{u}_N(v)dv=I_\vartheta \Phi^T(t)\Xi,\nonumber\\
u(t)=I_\vartheta\Phi^T(t)\Xi+u_0,\\
\label{Eq7}
I_\vartheta u=I_\vartheta^2 \Phi^T(t)\Xi+u_0t.
\label{Eq6}
\end{eqnarray}
Now, to obtain $\{\xi_i\}_{i=1}^N$ we define the residual functions by substituting Eqs. (\ref{Eq4})-(\ref{Eq6}) in Eq. (\ref{Eq3})
\begin{equation}
\hat{R}(t)=\kappa \Phi^T(t) \Xi - (I_\vartheta \Phi^T(t)\Xi +u_0)(1-I_\vartheta^T (t)\Xi -u_0 - I_\vartheta^2 \Phi^T(t)\Xi - u_0t).
\label{Eq8}
\end{equation}
The set of equations for obtaining the coefficients $\{\xi_i\}_{i=1}^N$ come from equalizing Eq. (\ref{Eq7}) to zero at $N$ interpolate nodes $\{t_i\}_{j=1}^N$ from $t_j=L(\frac{j}{N})^\rho,~~j=1,2,\cdots,N$
where $L$ is a last interpolate node and $\rho$ is a arbitrary parameter.
\begin{equation}
\hat{R}(t_j)=0,~~~~j=1,2,\cdots,N.
\label{Eq9}
\end{equation}

\section{Application}
\label{sec:5}
\begin{figure}
\centering \includegraphics[scale=0.5]{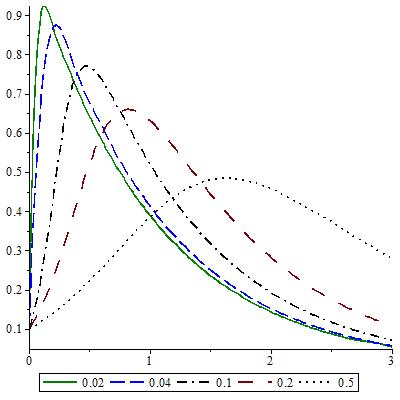}
\caption{Plot of approximate solutions of Eq. (\ref{Eq3}) for $u_0=0.1$ and $\kappa=0.02,0.04,0.1,0.2,0.5$. }
\label{Fig}
\end{figure}

We applied the method presented in this paper to examine the mathematical
structure of $u(t)$. Table (\ref{Tab1}) shows the maximum of $u(x)$ for some $\kappa$ and $u_0=0.1$ by using in comparison with exact solution and ADM solution by Wazwaz \cite{Wazwaz}. The resulting graph of Eq. (\ref{Eq3}) is shown in Fig. (\ref{Fig}).

\begin{table}[htbpH]
\caption{\scriptsize
A comparison of ADM\cite{Wazwaz} and the present method with exact values for $u_{max}$.
}
\scriptsize \begin{tabular*}{\columnwidth}{@{\extracolsep{\fill}}*{7}{|c|}}
\hline
$\kappa$&$u_{max}$ &$r_\omega$&$\rho$&$N$&$ICSRBF$&$ADM$\\
\hline
$0.02$&$0.9234271$ &$1$&$1.766000$&$15$&$0.92342716$&$0.9234270$\\
\hline
$0.04$&$0.8737199$ &$1$&$1.780000$&$18$&$0.8737193$&$0.8612401$\\
\hline
$0.1$&$0.7697414$ &$1$&$1.811000$&$18$&$0.7697414$&$0.7651130$\\
\hline
$0.2$&$0.6590503$ &$2$&$1.032770$&$18$&$0.6590493$&$0.6579123$\\
\hline
$0.5$&$0.4851902$ &$2$&$1.114035$&$27$&$0.4851903$&$0.4852823$\\
\hline
\end{tabular*}
\label{Tab1}
\end{table}

\begin{table}
\caption{Minimum value of $\|Res\|^2$ which is obtained with $r_\omega$ and $\rho$ for ICSRBF.}
\centering \begin{tabular}{|c|c|}
\hline 
$\kappa$&$\|Res\|^2$-ICSRBF\\
\hline
$0.5$&$2.11e-08$\\
\hline
$0.2$&$2.87e-07$\\
\hline
$0.1$&$1.34e-07$\\
\hline
$0.04$&$7.86e-05$\\
\hline
$0.02$&$4.41e-05$\\
\hline
\end{tabular} 
\label{Tab2}
\end{table}

The local support domain $r_\omega$ and arbitrary parameter $\rho$ must be specified by the user.An important
unsolved problem is to find a approach to determine the optimal size of $r_\omega$\cite{Wong}. The accuracy of these CSRBF depends on the choice of $r_\omega$ and $\rho$. By the meaning of residual function in case of Eq. (\ref{Eq7}), we try to minimize $\|R(t)\|^2$ by choosing good $r_\omega$ and $\rho$ parameters. We define $\|R(t)\|^2$ as 
\begin{equation}
\|R(t)\|^2=\int_0^b R^2(t)dt \simeq \sum_{j=0}^m \omega_j R^2 (\frac{L}{2}s_j+\frac{L}{2}),
\end{equation}
were 
\begin{equation}
\omega_j=\frac{L}{(1-s_j^2)(\frac{d}{dx}P_{m+1}(s)|_{s=s_j})^2},~~~~j=0,1,\cdots,m,\nonumber\\
P_{m+1}(s_j)=0,~~~~j=0,1,\cdots,m,\nonumber
\end{equation}
$P_{m+1}(x)$ is $(m+1)$th-order Legendre polynomial. Table (\ref{Tab2}) show the minimum of $\|R(t)\|^2$  which is obtained with local support domain $r_\omega$ and arbitrary parameter $\rho$ .

\section{ Conclusion}
\label{sec:6}
 A method has been presented for solving Volterra’s Population Model which is an integro-ordinary differential equation,
based on the compactly supported radial basis functions approximation. In this work, we applied two common  ICSRBF methods on the Volterra’s Population Model without converting it to an ordinary differential equation. We used 
$Wendland_{3,5}$ function. This function are proposed to provide an effective
but simple way to improve the convergence of the solution by collocation method. As appeared from the Figures, we have shown that,
when the constant $\kappa=\frac{c}{ab}$ is small, this type of population is relatively insensitive to toxins, and when $c=ab$ is large, population
of this type are extremely sensitive to toxins. Additionally, through the comparison with ADM, we have
showed that the  ICSRBF approach have good reliability and efficiency.



\section{\bf References}\label{refs}



\begin{thebibliography}{20}

\bibitem{Scudo}
F. Scudo ,
\newblock Vito Volterra and theoretical ecology,
\newblock {\em Theor Popul Biol}  
\textbf{2} (1971), 1--23.

\bibitem{TeBeest}
K. TeBeest,
\newblock Numerical and analytical solutions of Volterra’s population model,
\newblock {\em SIAM Rev}  
\textbf{39} (1997),484-493.

 \bibitem{Wazwaz}
A. M. Wazwaz,
\newblock {\em Analytical approximations and Padé approximants for Volterra’s population model},
\newblock {\em Appl Math Comput}  
\textbf{100} (1999),13--25.

\bibitem{DH9}
 B.~J. Noye, M. Dehghan, 
 \newblock  New explicit finite difference schemes for two-dimensional diffusion subject to specification of mass. \newblock {\em Numer Meth Par Diff Eq} \textbf{15} (1999) 521--534.

\bibitem{Small}
R. Small,
\newblock Population growth in a closed system,
\newblock {\em SIAM Rev}  \textbf{25} (1983) 93--95.

\bibitem{Al-Khaled}
K. Al-Khaled,
\newblock Numerical approximations for population growth models,
\newblock {\em Appl Math Comput}  \textbf{160} (2005) 865--873.

\bibitem{Parand1}
K. Parand , A. Rezaei, A. Taghavi,
\newblock Numerical approximations for population growth model by Rational Chebyshev and Hermite functions collocation approach: a comparison.,
\newblock {\em Math Methods Appl Sci}  \textbf{33} (2010) 2076--2086.

  \bibitem{DH7}
K. Parand, J.~A. Rad,
\newblock Kansa method for the solution of a parabolic equation with an unknown spacewise-dependent coefficient subject to an extra measurement. 
\newblock {\em Comp Phys Commun} \textbf{184} (2013) 582--595.
  
\bibitem{Ramezani}
M. Ramezani, M. Razzaghi, M. Dehghan,
\newblock Composite spectral functions for solving Volterra’s population model,
\newblock {\em Chaos Soliton Fract}  \textbf{34} (2007) 588--593.


   \bibitem{DH1}
M. Dehghan, A. Shokri, 
\newblock A numerical method for solution of the two-dimensional sine-gordon equation using the radial basis functions. 
\newblock {\em Math Comput Simul} \textbf{79} (2008) 700--715.

  \bibitem{DH6}
 M. Dehghan,A. Shokri, 
 \newblock A meshless method for numerical  solution of the one-dimensional wave equation with an integral condition  using radial basis functions. 
 \newblock {\em Numer Algorithms} \textbf{52} (2009) 461--477.

\bibitem{DH2}
M. Dehghan,A. Shokri,
 \newblock Numerical solution of the nonlinear klein–gordon equation using radial basis functions. 
\newblock {\em J Comput Appl Math} \textbf{230} (2009) 400--410.
  
  \bibitem{DH10}
 W. Bu,Y. Ting,Y.  Wu ,J. Yang,  
 \newblock Finite difference/finite element
  method for two-dimensional space and time fractional bloch–torrey
  equations. 
  \newblock {\em J Comput Phys} \textbf{293} (2015) 264--279.
  
  \bibitem{DH11}
H.~J. Choi, J.~R. Kweon, 
\newblock A finite element method for singular solutions
  of the navier–stokes equations on a non-convex polygon. 
  \newblock {\em J Comput Appl Math }\textbf{292} (2016) 342--362.


\bibitem{Parand6}
K. Parand, S. Abbasbandy, S. Kazem, J. A. Rad,
\newblock A novel application of radial basis functions for solving a model
of first-order integro-ordinary differential equation,
\newblock {\em Commun Nonlinear Sci Numer Simulat}  
\textbf{16} (2011) 4250--4258.

   \bibitem{DH4}
 J.~A. Rad, S. Kazem.,K. Parand, 
 \newblock A numerical solution of the nonlinear controlled duffing oscillator by radial basis function.
 \newblock {\em Comput Math Appl}
  \textbf{64} (2012) 2049--2065.
  
  \bibitem{DH5}
J.~A. Rad,S. Kazem,M. Shaban, K. Parand,A. Yildirim, 
\newblock Numerical solution of fractional differential equations with a tau method based on
  legendre and bernstein polynomials.
  \newblock {\em Math Meth Appl Sci}
  \textbf{37} (2014) 329--342.

\bibitem{Marzban}
HR. Marzban, S. Hoseini, M. Razzaghi,
\newblock Solution of Volterra’s population model via block-pulse functions and Lagrange-interpolating polynomials,
\newblock {\em Math
Methods Appl Sci}  
\textbf{32} (2009) 127--134.

\bibitem{Parand7}
K. Parand, Z. Delafkar, N. Pakniat, MK. Haji
\newblock Collocation method using Sinc and Rational Legendre functions for solving Volterra’s population
model,
\newblock {\em Commun Nonlinear Sci Numer Simul}  
\textbf{16} (2011) 1811--1819.

\bibitem{Momani}
S. Momani, R. Qaralleh, N. Pakniat, MK. Haji,
\newblock Numerical approximations and Padé approximants for a fractional population growth model,
\newblock {\em Appl Math Model}  
\textbf{31} (2007) 1907--1914.

\bibitem{Xu}
H. Xu,
\newblock Analytical approximations for a population growth model with fractional order,
\newblock {\em Commun Nonlinear Sci Numer Simul}  
\textbf{14} (2009) 1978--1983.

\bibitem{Wong}
S. M. Wong, Y. C. Hon, M. A. Golberg ,
\newblock Compactly supported radial basis function for shallow water equations,
\newblock {\em Appl  Math  Comput}  
\textbf{127} (2002) 79--101.

\bibitem{Wendland}
H. Wendland,
\newblock Piecewise polynomial, positive definite and compactly supported radial functions of minimal degree,
\newblock {\em Adv  Comput  Math}  
\textbf{4} (1995) 389–396.

\bibitem{MeshfreeFasshhauer}
G. E. Fasshauer,
\newblock Meshfree Approximation Methods With Matlab,
\newblock {\em World Scientific Publishing Co. (1995) Pte, Ltd}  
\textbf{4}  389–396.

\bibitem{Fasshauer}
G. E. Fasshauer,
\newblock On smoothing for multilevel approximation with radial basis functions, An approximation theory
IX, Vol. II: Coputational Aspects, CharlesK. Chui and L. L. Schumakher.,
\newblock {\em Vanderbilt University Press. (1999) Pte, Ltd}  
\textbf{4}  389–396.

 \bibitem{DH3}
A. Shokri, M. Dehghan, 
\newblock A not-a-knot meshless method using radial basis
  functions and predictor–corrector scheme to the numerical solution of
  improved boussinesq equation. 
  \newblock {\em Comput Phys Commun} \textbf{181} (2010) 1990--2000.
  
   \bibitem{DH8}
K. Rashidi,H. Adibi,J.~A. Rad, K. Parand, 
\newblock Application of meshfree
  methods for solving the inverse one-dimensional stefan problem,.
  \newblock {\em Eng Anal Bound Elem} \textbf{40} (2014) 1--21.

\end{thebibliography}
\end{document}